# A SET OF SEQUENCES IN NUMBER THEORY

by Florentin Smarandache
University of New Mexico
Gallup, NM 87301, USA

**Abstract**: New sequences are introduced in number theory, and for each one a general question: how many primes each sequence has.

**Keywords**: sequence, symmetry, consecutive, prime, representation of numbers.

**1991 MSC**: 11A67

**Introduction**.
74 new integer sequences are defined below, followed by references and some open questions.

## 1. Consecutive sequence:

1,12,123,1234,12345,123456,1234567,12345678,123456789,

12345678910,1234567891011,123456789101112,

12345678910111213,...

How many primes are there among these numbers?

In a general form, the Consecutive Sequence is considered in an arbitrary numeration base B.

Reference:

a) Student Conference, University of Craiova, Department of

   Mathematics, April 1979,  "Some problems in number

   theory" by Florentin Smarandache.

## 2. Circular sequence:

```
1,12,21,123,231,312,1234,2341,3412,4123,12345,23451,34512,45123,51234,
| |   | |         | |                 | |                            |
 ---   ---------   -----------------   ----------------------------
1   2      3              4                       5

123456,234561,345612,456123,561234,612345,1234567,2345671,3456712,...
|                                         | |
 ----------------------------------------   ---------------------  ...
                    6                                 7
```

How many primes are there among these numbers?

## 3. Symmetric sequence:

1,11,121,1221,12321,123321,1234321,12344321,123454321,

1234554321,12345654321,123456654321,1234567654321,

12345677654321,123456787654321,1234567887654321,

12345678987654321,123456789987654321,12345678910987654321,

1234567891010987654321,123456789101110987654321,

12345678910111110987654321,...

How many primes are there among these numbers?

In a general form, the Symmetric Sequence is considered in an arbitrary numeration base B.

References:

a) Arizona State University, Hayden Library, "The Florentin
   Smarandache papers" special collection, Tempe, AZ 85287-
   1006, USA.

b) Student Conference, University of Craiova, Department of
   Mathematics, April 1979, "Some problems in number
   theory" by Florentin Smarandache.

## 4. Deconstructive sequence:

```
1,23,456,7891,23456,789123,4567891,23456789,123456789,1234567891, ...
|       ||         ||       ||       | |      | |        ||
 ---------  --------- --------  --------  --------  -------  - ...
```

How many primes are there among these numbers?

Reference:

a) Arizona State University, Hayden Library, "The Florentin Smarandache papers" special collection, Tempe, AZ 85287-1006, USA.

## 5. Mirror sequence:

1,212,32123,4321234,543212345,65432123456,7654321234567,
876543212345678,98765432123456789,109876543212345678910,
1110987654321234567891011,...

Question: How many of them are primes?

## 6. Permutation sequence:

12, 1342, 135642, 13578642, 13579108642, 135791112108642,

1357911131412108642, 13579111315161412108642,

135791113151718161412108642, 13579111315171920181614121086 42, ...

Question: Is there any perfect power among these numbers?

(Their last digit should be:

either 2 for exponents of the form 4k+1,

either 8 for exponents of the form 4k+3, where k >= 0 .)

Conjecture: no!

## 7. Generalized permutation sequence:

If g(n), as a function, gives the number of digits of a(n),

and F if a permutation of g(n) elements, then:

$$a(n) = \overline{F(1)F(2)...F(g(n))} \ .$$

**8. Simple numbers**:

2,3,4,5,6,7,8,9,10,11,13,14,15,17,19,21,22,23,25,26,27,29,31,

33,34,35,37,38,39,41,43,45,46,47,49,51,53,55,57,58,61,62,65,

67,69,71,73,74,77,78,79,82,83,85,86,87,89,91,93,94,95,97,101,103,...

(A number n is called <simple number> if the product of its

proper divisors is less than or equal to n.)

Generally speaking, n has the form:

   n = p, or p^2, or p^3, or pq,  where p and q are distinct

primes.

How many primes, squares, and cubes are in this sequence?  What
interesting properties has it?

Reference:

a) Student Conference, University of Craiova, Department of

   Mathematics, April 1979,  "Some problems in number

   theory" by Florentin Smarandache.

## 9. Digital sum:

```
0,1,2,3,4,5,6,7,8,9,1,2,3,4,5,6,7,8,9,10,2,3,4,5,6,7,8,9,10,11,
|             | |                 | |                     |
 -----------------   ------------------   -------------------

3,4,5,6,7,8,9,10,11,12,4,5,6,7,8,9,10,11,12,13,5,6,7,8,9,10,11,12,13,14,...
|                  | |                     | |                        |
 ------------------   --------------------   ---------------------
```

(d$_s$(n) is the sum of digits.)

How many primes, squares, and cubes are in this sequence? What interesting properties has it?

## 10. Digital products:

```
0,1,2,3,4,5,6,7,8,9,0,1,2,3,4,5,6,7,8,9,0,2,4,6,8,19,12,14,16,18,
|                 | |                 | |                      |
 -----------------   -----------------   ----------------------

0,3,6,9,12,15,18,21,24,27,0,4,8,12,16,20,24,28,32,36,0,5,10,15,20,25,...
|                       | |                       | |
 -----------------------   -----------------------   --------------  ...
```

(d_p(n) is the product of digits.)

How many primes, squares, and cubes are in this sequence? What interesting properties has it?

## 11. Code puzzle:

151405,202315,2008180505,06152118,06092205,190924,1905220514,

0509070820,14091405,200514,051205220514,...

Using the following letter-to-number code:

```
 A  B  C  D  E  F  G  H  I  J  K  L  M  N  O  P  Q  R  S  T  U  V  W  X  Y  Z
01 02 03 04 05 06 07 08 09 10 11 12 13 14 15 16 17 18 19 20 21 22 23 24 25 26
```

then $c_p(n)$ = the numerical code of the spelling of n in English language; for exemple: 1 = ONE = 151405, etc.

Find a better codification (one sign only for each letter).

## 12. Pierced chain:

101,1010101,10101010101,101010101010101,1010101010101010101,

10101010101010101010101,101010101010101010101010101,...

```
(c(n) = 101 * 1 0001 0001 ... 0001 , for n >= 1)
              |  |  |  | ... |  |
              ---- ----      ----
               1    2         n-1
```

How many c(n)/101 are primes ?

References:

a) Arizona State University, Hayden Library, "The Florentin
   Smarandache papers" special collection, Tempe, AZ 85287-
   1006, USA.

b) Student Conference, University of Craiova, Department of
   Mathematics, April 1979, "Some problems in number
   theory" by Florentin Smarandache.

## 13. Divisor products:

1,2,3,8,5,36,7,64,27,100,11,1728,13,196,225,1024,17,5832,19,

8000,441,484,23,331776,125,676,729,21952,29,810000,31,32768,

1089,1156,1225,10077696,37,1444,1521,2560000,41,...

($P_d(n)$ is the product of all positive divisors of n.)

How many primes, squares, and cubes are in this sequence? What interesting properties has it?

## 14. Proper divisor products:

1,1,1,2,1,6,1,8,3,10,1,144,1,14,15,64,1,324,1,400,21,22,1,

13824,5,26,27,784,1,27000,1,1024,33,34,35,279936,1,38,39,

64000,1,...

($p_d(n)$ is the product of all positive divisors of n but n.)

How many primes, squares, and cubes are in this sequence? What interesting properties has it?

## 15. Square complements:

1,2,3,1,5,6,7,2,1,10,11,3,14,15,1,17,2,19,5,21,22,23,6,1,26,

3,7,29,30,31,2,33,34,35,1,37,38,39,10,41,42,43,11,5,46,47,3,

1,2,51,13,53,6,55,14,57,58,59,15,61,62,7,1,65,66,67,17,69,70,71,2,...

Definition:

 for each integer n to find the smallest integer k such that

 nk is a perfect square..

 (All these numbers are square free.)

    How many primes, squares, and cubes are in this sequence?  What interesting properties has it?

## 16. Cubic complements:

1,4,9,2,25,36,49,1,3,100,121,18,169,196,225,4,289,12,361,50,
441,484,529,9,5,676,1,841,900,961,2,1089,1156,1225,6,1369,
1444,1521,25,1681,1764,1849,242,75,2116,2209,36,7,20,...

Definition:

 for each integer n to find the smallest integer k such that

 nk is a perfect cub.

 (All these numbers are cube free.)

How many primes, squares, and cubes are in this sequence?  What interesting properties has it?

## 17. m-power complements:

Definition:

 for each integer n to find the smallest integer k such that

 nk is a perfect m-power (m => 2).

 (All these numbers are m-power free.)

How many primes, squares, and cubes are in this sequence?  What interesting properties has it?

Reference:

a) Arizona State University, Hayden Library, "The Florentin
   Smarandache papers" special collection, Tempe, AZ 85287-
   1006, USA.

### 18. Cube free sieve:

2,3,4,5,6,7,9,10,11,12,13,14,15,17,18,19,20,21,22,23,25,26,
28,29,30,31,33,34,35,36,37,38,39,41,42,43,44,45,46,47,49,50,
51,52,53,55,57,58,59,60,61,62,63,65,66,67,68,69,70,71,73,...

Definition:  from the set of natural numbers (except 0 and 1):

 - take off all multiples of $2^3$ (i.e. 8, 16, 24, 32, 40, ...)

 - take off all multiples of $3^3$

 - take off all multiples of $5^3$

 ... and so on (take off all multiples of all cubic primes).

(One obtains all cube free numbers.)

How many primes, squares, and cubes are in this sequence?  What interesting properties has it?

### 19. m-power free sieve:

Definition:  from the set of natural numbers (except 0 and 1)

 take off all multiples of $2^m$, afterwards all multiples of $3^m$, ...

 and so on (take off all multiples of all m-power primes, m >= 2).

(One obtains all m-power free numbers.)

How many primes, squares, and cubes are in this sequence?  What interesting properties has it?

## 20. Irrational root sieve:

2,3,5,6,7,10,11,12,13,14,15,17,18,19,20,21,22,23,24,26,28,
29,30,31,33,34,35,37,38,39,40,41,42,43,44,45,46,47,48,50,
51,52,53,54,55,56,57,58,59,60,61,62,63,65,66,67,68,69,70,71,72,73,...

Definition: from the set of natural numbers (except 0 and 1):

- take off all powers of $2^k$, k >= 2, (i.e. 4, 8, 16, 32, 64, ...)
- take off all powers of $3^k$, k >= 2;
- take off all powers of $5^k$, k >= 2;
- take off all powers of $6^k$, k >= 2;
- take off all powers of $7^k$, k >= 2;
- take off all powers of $10^k$, k >= 2;

... and so on (take off all k-powers, k >= 2, of all square free numbers).

One gets all square free numbers by the following method (sieve):

from the set of natural numbers (except 0 and 1):

- take off all multiples of $2^2$ (i.e. 4, 8, 12, 16, 20, ...)
- take off all multiples of $3^2$
- take off all multiples of $5^2$

... and so on (take off all multiples of all square primes);

one obtains, therefore:

2,3,5,6,7,10,11,13,14,15,17,19,21,22,23,26,29,30,31,33,34,
35,37,38,39,41,42,43,46,47,51,53,55,57,58,59,61,62,65,66,
67,69,70,71,... ,

which are used for irrational root sieve.

(One obtains all natural numbers those m-th roots, for any m >= 2, are irrational.)

How many primes, squares, and cubes are in this sequence? What interesting properties has it?

### 21. Odd sieve:

7,13,19,23,25,31,33,37,43,47,49,53,55,61,63,67,73,75,79,83,

85,91,93,97,...

(All odd numbers that are not equal to the difference of two primes.)

A sieve is used to get this sequence:

- substract 2 from all prime numbers and obtain a temporary sequence;

- choose all odd numbers that do not belong to the temporary one.

How many primes, squares, and cubes are in this sequence? What interesting properties has it?

## 22. Binary sieve:

1,3,5,9,11,13,17,21,25,27,29,33,35,37,43,49,51,53,57,59,65,
67,69,73,75,77,81,85,89,91,97,101,107,109,113,115,117,121,
123,129,131,133,137,139,145,149,...

(Starting to count on the natural numbers set at any step from 1:

  - delete every 2-nd numbers

  - delete, from the remaining ones, every 4-th numbers

  ... and so on:  delete, from the remaining ones, every

  $(2^k)$-th numbers, k = 1, 2, 3, ... .)

 Conjectures:

 - there are an infinity of primes that belong to this sequence;

 - there are an infinity of numbers of this sequence which

 are not prime.

**23. Trinary sieve**:

1,2,4,5,7,8,10,11,14,16,17,19,20,22,23,25,28,29,31,32,34,35,
37,38,41,43,46,47,49,50,52,55,56,58,59,61,62,64,65,68,70,71,
73,74,76,77,79,82,83,85,86,88,91,92,95,97,98,100,101,103,104,
106,109,110,112,113,115,116,118,119,122,124,125,127,128,130,
131,133,137,139,142,143,145,146,149,...

(Starting to count on the natural numbers set at any step from 1:
  - delete every 3-rd numbers
  - delete, from the remaining ones, every 9-th numbers
  ... and so on:  delete, from the remaining ones, every
  (3^k)-th numbers, k = 1, 2, 3, ... .)

 Conjectures:
 - there are an infinity of primes that belong to this sequence;
 - there are an infinity of numbers of this sequence which
 are not prime.

**24. n-ary sieve  (generalization, n >= 2):**

(Starting to count on the natural numbers set at any step from 1:
  - delete every n-th numbers;
  - delete, from the remaining ones, every (n^2)-th numbers;
  ... and so on:  delete, from the remaining ones, every
  (n^k)-th numbers, k = 1, 2, 3, ... .)

 Conjectures:
 - there are an infinity of primes that belong to this sequence;
 - there are an infinity of numbers of this sequence which

are not prime.

## 25. Consecutive sieve:

1,3,5,9,11,17,21,29,33,41,47,57,59,77,81,101,107,117,131,149.

153,173,191,209,213,239,257,273,281,321,329,359,371,401,417,

441,435,491,...

(From the natural numbers set:

- keep the first number,

  delete one number out of 2 from all remaining numbers;

- keep the first remaining number,

  delete one number out of 3 from the next remaining numbers;

- keep the first remaining number,

  delete one number out of 4 from the next remaining numbers;

... and so on, for step k (k >= 2):

- keep the first remaining number,

  delete one number out of k from the next remaining numbers;

... .)

This sequence is much less dense than the prime number sequence,

and their ratio tends to $p_n : n$ as n tends to infinity.

For this sequence we chosen to keep the first remaining number at all steps,

but in a more general case:

the kept number may be any among the remaining k-plet

(even at random).

How many primes, squares, and cubes are in this sequence?  What interesting properties has it?

## 26. General-sequence sieve:

Let $u_i > 1$, for $i = 1, 2, 3, \ldots$, a stricly increasing
positive integer sequence.  Then:

From the natural numbers set:

 - keep one number among 1, 2, 3, ..., $u_1 - 1$,
   and delete every $u_1$-th numbers;
 - keep one number among the next $u_2 - 1$ remaining numbers,
   and delete every $u_2$-th numbers;
 ... and so on, for step k (k >= 1):

 - keep one number among the next $u_k - 1$ remaining numbers,
   and delete every $u_k$-th numbers;
 ... .

   Problem:  study the relationship between sequence $u_i$, $i = 1, 2, 3, \ldots$,
   and the remaining sequence resulted from the general

   sieve.

   $u_i$, previously defined, is called sieve generator.

    How many primes, squares, and cubes are in this sequence?  What
interesting properties has it?

## 27. Digital sequences:

(This a particular case of sequences of sequences.)

General definition:

in any numeration base B, for any given infinite integer or rational sequence $S_1$, $S_2$, $S_3$, ..., and any digit D from 0 to B-1, it's built up a new integer sequence witch

   associates to $S_1$ the number of digits D of $S_1$ in base B,
   to $S_2$ the number of digits D of $S_2$ in base B, and so on...

For exemple, considering the prime number sequence in base 10, then the number of digits 1 (for exemple) of each prime number following their order is: 0,0,0,0,2,1,1,1,0,0,1,0,...

(the digit-1 prime sequence).

Second exemple if we consider the factorial sequence n! in base 10, then the number of digits 0 of each factorial number following their order is: 0,0,0,0,0,1,1,2,2,1,3,...

(the digit-0 factorial sequence).

Third exemple if we consider the sequence n^n in base 10, n=1,2,..., then the number of digits 5 of each term 1^1, 2^2, 3^3,..., following their order is: 0,0,0,1,1,1,1,0,0,0,...

(The digit-5 n^n sequence)

References:

a) E. Grosswald, University of Pennsylvania, Philadelphia,
   Letter to the Author, August 3, 1985;

b) R. K. Guy, University of Calgary, Alberta, Canada,

Letter to the Author, November 15, 1985;

c) Arizona State University, Hayden Library, "The Florentin Smarandache papers" special collection, Tempe, AZ 85287-1006, USA.

## 28. Construction sequences:

(This a particular case of sequences of sequences.)

General definition:

in any numeration base B, for any given infinite integer or rational sequence $S_1, S_2, S_3, ...,$ and any digits $D_1, D_2, ..., D_k$ (k < B), it's built up a new integer sequence such that each of its terms $Q_1 < Q_2 < Q_3 < ...$ is formed by these digits $D_1, D_2, ..., D_k$ only (all these digits are used), and matches a term $S_i$ of the previous sequence.

For exemple, considering in base 10 the prime number sequence, and the digits 1 and 7 (for exemple), we construct a written-only-with-these-digits (all these digits are used) prime number new sequence: 17,71,... (the digit-1-7-only prime sequence).

Second exemple, considering in base 10 the multiple of 3 sequence, and the digits 0 and 1, we construct a written-only-with-these-digits (all these digits are used) multiple of 3 new sequence: 1011,1101, 1110,10011,10101,10110,11001,11010,11100,... (the digit-0-1-only multiple of 3 sequence).

How many primes, squares, and cubes are in this sequence? What interesting properties has it?

References:

a) E. Grosswald, University of Pennsylvania, Philadelphia,

Letter to F. Smarandache, August 3, 1985;

b) R. K. Guy, University of Calgary, Alberta, Canada, Letter to F. Smarandache, November 15, 1985;

c) Arizona State University, Hayden Library, "The Florentin Smarandache papers" special collection, Tempe, AZ 85287-1006, USA.

## 29. General residual sequence:

$(x + C_1)...(x + C_{F(m)})$, m = 2, 3, 4, ...,

where $C_i$, $1 \le i \le F(m)$, forms a reduced set of residues mod m,

x is an integer, and F is Euler's totient.

The General Residual Sequence is induced from the

The Residual Function (see <Libertas Mathematica>):

Let $L : Z \times Z \rightarrow Z$ be a function defined by

$L(x,m) = (x + C_1)...(x + C_{F(m)})$,

where $C_i$, $1 \le i \le F(m)$, forms a reduced set of residues mod m,

$m \ge 2$, x is an integer, and F is Euler's totient.

The Residual Function is important because it generalizes the classical theorems by Wilson, Fermat, Euler, Wilson,

Gauss, Lagrange, Leibnitz, Moser, and Sierpinski all together.

For x=0 it's obtained the following sequence:

$L(m) = C_1 ... C_{F(m)}$, where m = 2, 3, 4, ...

(the product of all residues of a reduced set mod m):

1,2,3,24,5,720,105,2240,189,3628800,385,479001600,19305,

896896,2027025,20922789888000,85085,6402373705728000,

8729721,47297536000,1249937325,...

which is found in "The Handbook of Integer Sequences", by N.

J. A. Sloane, Academic Press, USA, 1973.

The Residual Function extends it.

How many primes, squares, and cubes are in this sequence? What interesting properties has it?

## 30. (Inferior) prime part:

2,3,3,5,5,7,7,7,7,11,11,13,13,13,13,17,17,19,19,19,19,23,23,

23,23,23,23,29,29,31,31,31,31,31,31,37,37,37,37,41,41,43,43,

43,43,47,47,47,47,47,47,53,53,53,53,53,53,59, ...

(For any positive real number n one defines $p_p(n)$ as the
largest prime number less than or equal to n.)

## 31. (Superior) prime part:

2,2,2,3,5,5,7,7,11,11,11,11,13,13,17,17,17,17,19,19,23,23,

23,23,29,29,29,29,29,29,31,31,37,37,37,37,37,37,41,41,41,

41,43,43,47,47,47,47,53,53,53,53,53,53,59,59,59,59,59,59,61,...

(For any positive real number n one defines $P_p(n)$ as the
smallest prime number greater than or equal to n.)

Study these sequences.

Reference:

a) Arizona State University, Hayden Library, "The Florentin
   Smarandache papers" special collection, Tempe, AZ 85287-
   1006, USA.

## 32. (Inferior) square part:

0,1,1,1,4,4,4,4,4,9,9,9,9,9,9,9,16,16,16,16,16,16,16,16,16,

25,25,25,25,25,25,25,25,25,25,25,36,36,36,36,36,36,36,36,

36,36,36,36,49,49,49,49,49,49,49,49,49,49,49,49,49,49,64,64,...

(The largest square less than or equal to n.)

## 33. (Superior) square part:

0,1,4,4,4,9,9,9,9,9,16,16,16,16,16,16,16,25,25,25,25,25,25,

25,25,25,36,36,36,36,36,36,36,36,36,36,36,49,49,49,49,49,49,

49,49,49,49,49,49,49,64,64,64,64,64,64,64,64,64,64,64,64,64,

64,64,81,81,...

(The smallest square greater than or equal to n.)

Study these sequences.

**34. (Inferior) cube part:**

0,1,1,1,1,1,1,1,8,8,8,8,8,8,8,8,8,8,8,8,8,8,8,8,8,8,27,27,

27,27,27,27,27,27,27,27,27,27,27,27,27,27,27,27,27,27,27,

27,27,27,27,27,27,27,27,27,27,27,27,27,27,64,64,64,...

(The largest cube less than or equal to n.)

**35. (Superior) cube part:**

0,1,8,8,8,8,8,8,8,27,27,27,27,27,27,27,27,27,27,27,27,27,

27,27,27,27,27,27,64,64,64,64,64,64,64,64,64,64,64,64,64,

64,64,64,64,64,64,64,64,64,64,64,64,64,64,64,64,64,64,64,

64,64,64,64,64,125,125,125,..

(The smalest cube greater than or equal to n.)

Study these sequences.

## 36. (Inferior) factorial part:

1,2,2,2,2,6,6,6,6,6,6,6,6,6,6,6,6,6,6,6,6,6,24,24,24,24,
24,24,24,24,24,24,24,24,24,24,24,24,24,24,...

(F$_p$(n) is the largest factorial less than or equal to n.)

## 37. (Superior) factorial part:

1,2,6,6,6,6,24,24,24,24,24,24,24,24,24,24,24,24,24,24,24,24,
24,24,120,120,120,120,120,120,120,120,120,120,120,...

(f$_p$(n) is the smallest factorial greater than or equal to n.)

Study these sequences.

## 38. Double factorial complements:

1,1,1,2,3,8,15,1,105,192,945,4,10395,46080,1,3,2027025,2560,
34459425,192,5,3715891200,13749310575,2,81081,1961990553600,
35,23040,213458046676875,128,6190283353629375,12,...

(For each n to find the smallest k such that nk is a double
 factorial, i.e. nk = either 1*3*5*7*9*...*n if n is odd,
                    either 2*4*6*8*...*n if n is even.)

Study this sequence in interrelation with Smarandache
function { S(n) is the smallest integer such that S(n)! is
divisible by n ).

## 39. Prime additive complements:

1,0,0,1,0,1,0,3,2,1,0,1,0,3,2,1,0,1,0,3,2,1,0,5,4,3,2,1,0,
1,0,5,4,3,2,1,0,3,2,1,0,1,0,3,2,1,0,5,4,3,2,1,0,...

(For each n to find the smallest k such that n+k is prime.)

Remark: is it possible to get as large as we want
but finite decreasing k, k-1, k-2, ..., 2, 1, 0 (odd k)
sequence included in the previous sequence -- i.e. for any
even integer are there two primes those difference is equal
to it? I conjecture the answer is negative.

Reference:

a) Arizona State University, Hayden Library, "The Florentin
   Smarandache papers" special collection, Tempe, AZ 85287-
   1006, USA.

**40. Prime base**:

0,1,10,100,101,1000,1001,10000,10001,10010,10100,100000,

100001,1000000,1000001,1000010,1000100,10000000,10000001,

100000000,100000001,100000010,100000100,1000000000,1000000001,

1000000010,1000000100,1000000101,...

(Each number n written in the prime base.)

(I define over the set of natural numbers the following infinite

base: $p_0 = 1$, and for $k \geq 1$ $p_k$ is the k-th prime number.)

Every positive integer A may be uniquely written in

the prime base as:

$$A = \overline{(a_n \ldots a_1 a_0)}_{(SP)} \stackrel{def}{===} \sum_{i=0}^{n} a_i p_i, \text{ with all } a_i = 0 \text{ or } 1, \text{ (of course } a_n = 1\text{)},$$

in the following way:

- if $p_n \leq A < p_{n+1}$ then $A = p_n + r_1$;
- if $p_m \leq r_1 < p_{m+1}$ then $r_1 = p_m + r_2$, $m < n$;

and so on untill one obtains a rest $r_j = 0$.

Therefore, any number may be written as a sum of prime numbers + e,

where e = 0 or 1.

If we note by p(A) the superior part of A (i.e. the largest

prime less than or equal to A), then

A is written in the prime base as:

   A = p(A) + p(A-p(A)) + p(A-p(A)-p(A-p(A))) + ... .

This base is important for partitions with primes.

    How many primes, squares, and cubes are in this sequence?   What
interesting properties has it?

## 41. Square base:

0,1,2,3,10,11,12,13,20,100,101,102,103,110,111,112,1000,1001,
1002,1003,1010,1011,1012,1013,1020,10000,10001,10002,10003,
10010,10011,10012,10013,10020,10100,10101,100000,100001,100002,
100003,100010,100011,100012,100013,100020,100100,100101,100102,
100103,100110,100111,100112,101000,101001,101002,101003,101010,
101011,101012,101013,101020,101100,101101,101102,1000000,...

(Each number n written in the square base.)

(I define over the set of natural numbers the following infinite
base: for $k \geq 0$  $s_k = k^2$.)

Every positive integer A may be uniquely written in

the square base as:

$$A = \overline{(a_n \ldots a_1 a_0)}_{(S2)} \stackrel{def}{===} \sum_{i=0}^{n} a_i s_i, \text{ with } a_i = 0 \text{ or } 1 \text{ for } i \geq 2,$$

$0 \leq a_0 \leq 3$, $0 \leq a_1 \leq 2$, and of course $a_n = 1$,

in the following way:

- if $s_n \leq A < s_{n+1}$ then $A = s_n + r_1$;
- if $s_m \leq r_1 < p_{m+1}$ then $r_1 = s_m + r_2$, $m < n$;

and so on untill one obtains a rest $r_j = 0$.

Therefore, any number may be written as a sum of squares

(1 not counted as a square -- being obvious) + e, where

e = 0, 1, or 3.

If we note by s(A) the superior square part of A (i.e. the

largest square less than or equal to A), then A is written

in the square base as:

$$A = s(A) + s(A-s(A)) + s(A-s(A)-s(A-s(A))) + \ldots .$$

This base is important for partitions with squares.

How many primes, squares, and cubes are in this sequence? What interesting properties has it?

## 42. m-power base  (generalization):

(Each number n written in the m-power base,

where m is an integer >= 2.)

(I define over the set of natural numbers the following infinite

m-power base:  for $k \geq 0$  $t_k = k^m$.)

Every positive integer A may be uniquely written in

the m-power base as:

$$A = \overline{(a_n \ldots a_1 a_0)}_{(SM)} \stackrel{def}{===} \sum_{i=0}^{n} a_i t_i \text{, with } a_i = 0 \text{ or } 1 \text{ for } i \geq m,$$

$0 \leq a_i \leq \lfloor ((i+2)^m - 1) / (i+1)^m \rfloor$ (integer part)

for $i = 0, 1, \ldots, m-1$, $a_i = 0$ or $1$ for $i \geq m$, and of course $a_n = 1$,
in the following way:
- if $t_n \leq A < t_{n+1}$  then $A = t_n + r_1$ ;
- if $t_m \leq r_1 < t_{m+1}$  then $r_1 = t_m + r_2$, $m < n$;
 and so on untill one obtains a rest $r_j = 0$.

Therefore, any number may be written as a sum of m-powers

(1 not counted as an m-power -- being obvious) + e, where

$e = 0, 1, 2, \ldots,$ or $2^m-1$.

If we note by t(A) the superior m-power part of A (i.e. the

largest m-power less than or equal to A), then A is written in the

m-power base as:

   $A = t(A) + t(A-t(A)) + t(A-t(A)-t(A-t(A))) + \ldots$

This base is important for partitions with m-powers.

   How many primes, squares, and cubes are in this sequence?  What
interesting properties has it?

### 43. Generalized base:

(Each number n written in the generalized base.)

(I define over the set of natural numbers the following infinite generalized base: $1 = g_0 < g_1 < \ldots < g_k < \ldots$ .)

Every positive integer A may be uniquely written in the generalized base as:

$$A = \overline{(a_n \ldots a_1 a_0)}_{(SG)} \stackrel{def}{===} \sum_{i=0}^{n} a_i g_i, \text{ with } 0 \leq a_i \leq \lfloor (g_{i+1} - 1) / g_i \rfloor$$

(integer part) for $i = 0, 1, \ldots, n$, and of course $a_n \geq 1$, in the following way:

- if $g_n \leq A < g_{n+1}$ then $A = g_n + r_1$ ;
- if $g_m \leq r_1 < g_{m+1}$ then $r_1 = g_m + r_2$, $m < n$;

and so on untill one obtains a rest $r_j = 0$.

If we note by $g(A)$ the superior generalized part of A (i.e. the largest $g_i$ less than or equal to A), then A is written in the m-power base as:

$$A = g(A) + g(A-g(A)) + g(A-g(A)-g(A-g(A))) + \ldots$$

This base is important for partitions: the generalized base may be any infinite integer set (primes, squares, cubes, any m-powers, Fibonacci/Lucas numbers, Bernoully numbers, etc.) those partitions are studied.

A particular case is when the base verifies: $2g_i \geq g_{i+1}$ for any i, and $g_0 = 1$, because all coefficients of a written number in this base will be 0 or 1.

Remark: another particular case: if one takes $g_i = p^{i-1}$, $i = 1, 2, 3, \ldots$, p an integer $\geq 2$, one gets the representation of a number in the numerical base p {p may be 10 (decimal), 2 (binar), 16 (hexadecimal),

etc.}.

How many primes, squares, and cubes are in this sequence?  What interesting properties has it?

## 44. Factorial quotients:

1,1,2,6,24,1,720,3,80,12,3628800,2,479001600,360,8,45,

20922789888000,40,6402373705728000,6,240,1814400,

1124000727777607680000,1,145152,239500800,13440,180,

304888344611713860501504000000,...

(For each n to find the smallest k such that nk is a factorial number.)

Study this sequence in interrelation with Smarandache function.

Reference:

a) Arizona State University, Hayden Library, "The Florentin
   Smarandache papers" special collection, Tempe, AZ 85287-
   1006, USA.

## 45. Double factorial numbers:

1,2,3,4,5,6,7,4,9,10,11,6,13,14,5,6,17,12,19,10,7,22,23,6,

15,26,9,14,29,10,31,8,11,34,7,12,37,38,13,10,41,14,43,22,9,

46,47,6,21,10,...

($d_f(n)$ is the smallest integer such that $d_f(n)!!$ is a multiple of n.)

Study this sequence in interrelation with Smarandache function.

Reference:

a) Arizona State University, Hayden Library, "The Florentin Smarandache papers" special collection, Tempe, AZ 85287-1006, USA.

## 46. Primitive numbers (of power 2):

2,4,4,6,8,8,8,10,12,12,14,16,16,16,16,18,20,20,22,24,24,24,

26,28,28,30,32,32,32,32,32,34,36,36,38,40,40,40,42,44,44,46,

48,48,48,48,50,52,52,54,56,56,56,58,60,60,62,64,64,64,64,64,64,66,...

($S_2(n)$ is the smallest integer such that $S_2(n)!$ is divisible by $2^n$.)

Curious property: this is the sequence of even numbers,

each number being repeated as many times as its exponent

(of power 2) is.

This is one of irreductible functions, noted $S_2(k)$, which helps to calculate the Smarandache function.

How many primes, squares, and cubes are in this sequence? What interesting properties has it?

## 47. Primitive numbers (of power 3):

3,6,9,9,12,15,18,18,21,24,27,27,27,30,33,36,36,39,42,45,45,

48,51,54,54,54,57,60,63,63,66,69,72,72,75,78,81,81,81,81,84,

87,90,90,93,96,99,99,102,105,108,108,108,111,...

($S_3(n)$ is the smallest integer such that $S_3(n)!$ is divisible by $3^n$.)

Curious property: this is the sequence of multiples of 3,

each number being repeated as many times as its exponent

(of power 3) is.

This is one of irreducible functions, noted $S_3(k)$, which helps to calculate the Smarandache function.

How many primes, squares, and cubes are in this sequence? What interesting properties has it?

## 48. Primitive numbers (of power p, p prime) {generalization}:

($S_p(n)$ is the smallest integer such that $S_p(n)!$ is divisible by $p^n$.)

Curious property: this is the sequence of multiples of p, each number being repeated as many times as its exponent (of power p) is.

These are the irreducible functions, noted $S_p(k)$, for any prime number p, which helps to calculate the Smarandache function.

How many primes, squares, and cubes are in this sequence? What interesting properties has it?

## 49. Square residues:

1,2,3,2,5,6,7,2,3,10,11,6,13,14,15,2,17,6,19,10,21,22,23,6,

5,26,3,14,29,30,31,2,33,34,35,6,37,38,39,10,41,42,43,22,15,

46,47,6,7,10,51,26,53,6,14,57,58,59,30,61,62,21,...

($s_r(n)$ is the largest square free number which divides n.)

Or, $s_r(n)$ is the number n released of its squares:
if $n = (p_1 \wedge a_1) * ... * (p_r \wedge a_r)$, with all $p_i$ primes and all $a_i \geq 1$,
then $s_r(n) = p_1 * ... * p_r$.

Remark: at least the $(2^2)*k$-th numbers (k = 1, 2, 3, ...)

are released of their squares;

and more general: all $(p^2)*k$-th numbers (for all p prime,

and k = 1, 2, 3, ...) are released of their squares.

How many primes, squares, and cubes are in this sequence? What interesting properties has it?

## 50. Cubical residues:

1,2,3,4,5,6,7,4,9,10,11,12,13,14,15,4,17,18,19,20,21,22,23,

12,25,26,9,28,29,30,31,4,33,34,35,36,37,38,39,20,41,42,43,

44,45,46,47,12,49,50,51,52,53,18,55,28,...

($c_r(n)$ is the largest cube free number which divides n.)

Or, $c_r(n)$ is the number n released of its cubicals:
if $n = (p_1 \char`\^ a_1) * ... * (p_r \char`\^ a_r)$, with all $p_i$ primes and all $a_i >= 1$,
then $c_r(n) = (p_1 \char`\^ b_1) * ... * (p_r \char`\^ b_r)$, with all $b_i = \min\{2, a_i\}$.

Remark: at least the (2^3)*k-th numbers (k = 1, 2, 3, ...)

are released of their cubicals;

and more general: all (p^3)*k-th numers (for all p prime,

and k = 1, 2, 3, ...) are released of their cubicals.

How many primes, squares, and cubes are in this sequence? What interesting properties has it?

## 51. m-power residues   (generalization):

$m_r(n)$ is the largest m-power free number which divides n.

Or, $m_r(n)$ is the number n released of its m-powers:
if $n = (p_1 \wedge a_1) * \ldots * (p_r \wedge a_r)$, with all $p_i$ primes and all $a_i \geq 1$,
then $m_r(n) = (p_1 \wedge b_1) * \ldots * (p_r \wedge b_r)$, with all $b_i = \min\{m-1, a_i\}$.

Remark:  at least the $(2^m)*k$-th numbers (k = 1, 2, 3, ...)

are released of their m-powers;

and more general:  all $(p^m)*k$-th numers (for all p prime,

and k = 1, 2, 3, ...) are released of their m-powers.

How many primes, squares, and cubes are in this sequence?  What interesting properties has it?

## 52. Exponents (of power 2):

0,1,0,2,0,1,0,3,0,1,0,2,0,1,0,4,0,1,0,2,0,1,0,2,0,1,0,2,0,

1,0,5,0,1,0,2,0,1,0,3,0,1,0,2,0,1,0,3,0,1,0,2,0,1,0,3,0,1,

0,2,0,1,0,6,0,1,...

($e_2(n)$ is the largest exponent (of power 2) which divides n .)

Or, $e_2(n) = k$ if $2^k$ divides n but $2^{(k+1)}$ does not.

## 53. Exponents (of power 3):

0,0,1,0,0,1,0,0,2,0,0,1,0,0,1,0,0,2,0,0,1,0,0,1,0,0,3,0,0,

1,0,0,1,0,0,2,0,0,1,0,0,1,0,0,2,0,0,1,0,0,1,0,0,2,0,0,1,0,

0,1,0,0,2,0,0,1,0,...

($e_3(n)$ is the largest exponent (of power 3) which divides n.)

Or, $e_3(n) = k$ if $3^k$ divides n but $3^{(k+1)}$ does not.

## 54. Exponents (of power p) {generalization}:

($e_p(n)$ is the largest exponent (of power p) which divides n, where p is an integer >= 2 .)

Or, $e_p(n) = k$ if $p^k$ divides n but $p^{(k+1)}$ does not.

Study these sequences.

Reference:

a) Arizona State University, Hayden Library, "The Florentin Smarandache papers" special collection, Tempe, AZ 85287-1006, USA.

## 55. Pseudo-primes:

2,3,5,7,11,13,14,16,17,19,20,23,29,30,31,32,34,35,37,38,41,

43,47,50,53,59,61,67,70,71,73,74,76,79,83,89,91,92,95,97,98,

101,103,104,106,107,109,110,112,113,115,118,119,121,124,125,

127,128,130,131,133,134,136,137,139,140,142,143,145,146, ...

(A number is pseudo-prime if some permutation of the digits

 is a prime number, including the identity permutation.)

(Of course, all primes are pseudo-primes,

 but not the reverse!)

How many primes, squares, and cubes are in this sequence? What interesting properties has it?

## 56. Pseudo-squares:

1,4,9,10,16,18,25,36,40,46,49,52,61,63,64,81,90,94,100,106,
108,112,121,136,144,148,160,163,169,180,184,196,205,211,225,
234,243,250,252,256,259,265,279,289,295,297,298,306,316,324,
342,360,361,400,406,409,414,418,423,432,441,448,460,478,481,
484,487,490,502,520,522,526,529,562,567,576,592,601,603,604,
610,613,619,625,630,631,640,652,657,667,675,676,691,729,748,
756,765,766,784,792,801,810,814,829,841,844,847,874,892,900,
904,916,925,927,928,940,952,961,972,982,1000, ...

(A number is a pseudo-square if some permutation of the

digits is aperfect square, including the identity permutation.)

(Of course, all perfect squares are pseudo-squares,

but not the reverse!)

One listed all pseudo-squares up to 1000.

How many primes, squares, and cubes are in this sequence? What interesting properties has it?

## 57. Pseudo-cubes:

1,8,10,27,46,64,72,80,100,125,126,152,162,207,215,216,251,
261,270,279,297,334,343,406,433,460,512,521,604,612,621,
640,702,720,729,792,800,927,972,1000,...

(A number is a pseudo-cube if some permutation of the digits
is a cube, including the identity permutation.)

(Of course, all perfect cubes are pseudo-cubes,
but not the reverse!)

One listed all pseudo-cubes up to 1000.

## 58. Pseudo-m-powers:

(A number is a pseudo-m-power if some permutation of the digits is an
m-power, including the identity permutation; m >= 2.)

Study these sequences.

Reference:
a) Arizona State University, Hayden Library, "The Florentin
   Smarandache papers" special collection, Tempe, AZ 85287-
   1006, USA.

### 59. Pseudo-factorials:

1,2,6,10,20,24,42,60,100,102,120,200,201,204,207,210,240,270,
402,420,600,702,720,1000,1002,1020,1200,2000,2001,2004,2007,
2010,2040,2070,2100,2400,2700,4002,4005,4020,4050,4200,4500,
5004,5040,5400,6000,7002,7020,7200,...

(A number is a pseudo-factorial if some permutation of the digits
is a factorial number, including the identity permutation.)

(Of course, all factorials are pseudo-factorials,
but not the reverse!)
One listed all pseudo-factorials up to 10000.

Procedure to obtain this sequence:
 - calculate all factorials with one digit only (1!=1, 2!=2,
   and 3!=6), this is line_1 (of one digit pseudo-factorials):
   1,2,6;
 - add 0 (zero) at the end of each element of line_1,
   calculate all factorials with two digits (4!=24 only)
   and all permutations of their digits:
   this is line_2 (of two digits pseudo-factorials):
   10,20,60; 24, 42;
 - add 0 (zero) at the end of each element of line_2 as well
   as anywhere in between their digits,
   calculate all factorials with three digits (5!=120,
   and 6!=720) and all permutations of their digits:
   this is line_3 (of three digits pseudo-factorials):
   100,200,600,240,420,204,402; 120,720, 102,210,201,702,270,720;
 and so on ...
 to get from line_k to line_(k+1) do:

- add 0 (zero) at the end of each element of line_k as well

  as anywhere in between their digits,

  calculate all factorials with (k+1) digits

  and all permutations of their digits;

The set will be formed by all line_1 to the last line elements in an increasing order.

How many primes, squares, and cubes are in this sequence?  What interesting properties has it?

## 60. Pseudo-divisors:

1,10,100,1,2,10,20,100,200,1,3,10,30,100,300,1,2,4,10,20,40,
100,200,400,1,5,10,50,100,500,1,2,3,6,10,20,30,60,100,200,
300,600,1,7,10,70,100,700,1,2,4,8,10,20,40,80,100,200,400,
800,1,3,9,10,30,90,100,300,900,1,2,5,10,20,50,100,200,500,1000,...

(The pseudo-divisors of n.)

(A number is a pseudo-divisor of n if some permutation of the
digits is a divisor of n, including the identity permutation.)

(Of course, all divisors are pseudo-divisors,
but not the reverse!)

A strange property:  any integer has an infinity of
pseudo-divisors !!
because 10...0 becomes 0...01 = 1, by a circular permutation
of its digits, and 1 divides any integer !

One listed all pseudo-divisors up to 1000 for the numbers 1, 2, 3,
..., 10.

Procedure to obtain this sequence:
 - calculate all divisors with one digit only,
   this is line_1 (of one digit pseudo-divisors);
 - add 0 (zero) at the end of each element of line_1,
   calculate all divisors with two digits
   and all permutations of their digits:
   this is line_2 (of two digits pseudo-divisors);
 - add 0 (zero) at the end of each element of line_2 as well
   as anywhere in between their digits,
   calculate all divisors with three digits

and all permutations of their digits:

   this is line_3 (of three digits pseudo-divisors);

  and so on ...

  to get from line_k to line_(k+1) do:

  - add 0 (zero) at the end of each element of line_k as well

    as anywhere in between their digits,

    calculate all divisors with (k+1) digits

    and all permutations of their digits;

  The set will be formed by all line_1 to the last line elements

  in an increasing order.

   How many primes, squares, and cubes are in this sequence?  What interesting properties has it?

## 61. Pseudo-odd numbers:

1,3,5,7,9,10,11,12,13,14,15,16,17,18,19,21,23,25,27,29,30,

31,32,33,34,35,36,37,38,39,41,43,45,47,49,50,51,52,53,54,

55,56,57,58,59,61,63,65,67,69,70,71,72,73,74,75,76,...

(Some permutation of digits is an odd number.)

How many primes, squares, and cubes are in this sequence? What interesting properties has it?

## 62. Pseudo-triangular numbers:

1,3,6,10,12,15,19,21,28,30,36,45,54,55,60,61,63,66,78,82,

87,91,...

(Some permutation of digits is a triangular number.)

A triangular number has the general form: n(n+1)/2.

How many primes, squares, and cubes are in this sequence? What interesting properties has it?

## 63. Pseudo-even numbers:

0,2,4,6,8,10,12,14,16,18,20,21,22,23,24,25,26,27,28,29,30,

32,34,36,38,40,41,42,43,44,45,46,47,48,49,50,52,54,56,58,

60,61,62,63,64,65,66,67,68,69,70,72,74,76,78,80,81,82,83,

84,85,86,87,88,89,90,92,94,96,98,100,...

(The pseudo-even numbers.)

(A number is a pseudo-even number if some permutation of

 the digits is a even number, including the identity permutation.)

(Of course, all even numbers are pseudo-even numbers,

 but not the reverse!)

A strange property:  an odd number can be a pseudo-even number!

One listed all pseudo-even numbers up to 100.

How many primes, squares, and cubes are in this sequence?  What interesting properties has it?

**64. Pseudo-multiples (of 5)**:

0,5,10,15,20,25,30,35,40,45,50,51,52,53,54,55,56,57,58,59,

60,65,70,75,80,85,90,95,100,101,102,103,104,105,106,107,

108,109,110,115,120,125,130,135,140,145,150,151,152,153,

154,155,156,157,158,159,160,165,...

(The pseudo-multiples of 5.)

(A number is a pseudo-multiple of 5 if some permutation of the digits is
 a multiple of 5, including the identity permutation.)

(Of course, all multiples of 5 are pseudo-multiples,
 but not the reverse!)

**65. Pseudo-multiples of p (p is an integer >= 2) {generalizations}:**

   (The pseudo-multiples of p.)

(A number is a pseudo-multiple of p if some permutation of the digits is
 a multiple of p, including the identity permutation.)

(Of course, all multiples of p are pseudo-multiples,
 but not the reverse!)

 Procedure to obtain this sequence:
  - calculate all multiples of p with one digit only (if any),
    this is line_1 (of one digit pseudo-multiples of p);
  - add 0 (zero) at the end of each element of line_1,
    calculate all multiples of p with two digits (if any)
    and all permutations of their digits:

this is line_2 (of two digits pseudo-multiples of p);

- add 0 (zero) at the end of each element of line_2 as well

  as anywhere in between their digits,

  calculate all multiples with three digits (if any)

  and all permutations of their digits:

  this is line_3 (of three digits pseudo-multiples of p);

and so on ...

to get from line_k to line_(k+1) do:

- add 0 (zero) at the end of each element of line_k as well

  as anywhere in between their digits,

  calculate all multiples with (k+1) digits (if any)

  and all permutations of their digits;

The set will be formed by all line_1 to the last line elements in an increasing order.

     How many primes, squares, and cubes are in this sequence?  What interesting properties has it?

     Reference:

a) Arizona State University, Hayden Library, "The Florentin

   Smarandache papers" special collection, Tempe, AZ 85287-

   1006, USA.

## 66. The Generalized Palindrome:

has one of the forms:

a(1)a(2)...a(n-1)a(n)a(n-1)...a(2)a(1) or

a(1)a(2)...a(n-1)a(n)a(n)a(n-1)...a(2)a(1),

where all a(k) are positive integers of one or more digits, and all above a(k) integers are concatenated.

(When all a(k) integers have a digit only, one gets the classical definition of the palindrome.)

Obviously, when n=1 in the first case, one can consider any positive integer as a SGP because, say 1743902 = (1743902), i.e. a(1) = 1743902. But let's take in the first case n > 1.

Examples of GP:

1457567145 because 1457567145 = (145)(7)(56)(7)(145),

also 145756567145 because 145756567145 = (145)(7)(56)(56)(7)(145).

Question: how many terms from the following sequences (known as Smarandache symmetric sequences) are primes?

a) 121, 12321, ..., 12...898...21, 12...9109...21, 12...91011109...21, etc.

b) 1221, 123321, 123...8998...21, 123...89101098...21, 123...891011111098...21, etc.

Reference:

a) F.Smarandache, Properties of numbers, University of Craiova, 1972.

**67. Constructive set (of digits 1,2):**

1,2,11,12,21,22,111,112,121,122,211,212,221,222,1111,1112,

1121,1122,1211,1212,1221,1222,21112112,2121,2122,2211,2212,

2221,2222,...

(Numbers formed by digits 1 and 2 only.)

Definition:

a1) 1, 2 belong to S;

a2) if a, b belong to S, then $\overline{ab}$ belongs to S too;

a3) only elements obtained by rules a1) and a2) applied a
    finite number of times belong to S.

Remark:

- there are 2^k numbers of k digits in the sequence, for
  k = 1, 2, 3, ... ;
- to obtain from the k-digits number group the (k+1)-digits
  number group, just put first the digit 1 and second the
  digit 2 in the front of all k-digits numbers.

```
Constructive set (of digits 1,2,3):

1,2,3,11,12,13,21,22,23,31,32,33,111,112,113,121,122,123,

131,132,133,211,212,213,221,222,223,231,232,233,311,312,

313,321,322,323,331,332,333,...

(Numbers formed by digits 1, 2, and 3 only.)

Definition:

a1) 1, 2, 3 belong to S;

a2) if a, b belong to S, then $\overline{ab}$ belongs to S too;

a3) only elements obtained by rules a1) and a2) applied

    a finite number of times belong to S.

Remark:

 - there are $3^k$ numbers of k digits in the sequence, for

   k = 1, 2, 3, ... ;

 - to obtain from the k-digits number group the (k+1)-digits

   number group, just put first the digit 1, second the digit 2,

   and third the digit 3 in the front of all k-digits numbers.
```

## 68. Generalizated constructive set:

(Numbers formed by digits $d_1, d_2, ..., d_m$ only,
 all $d_i$ being different each other, $1 \le m \le 9$.)

Definition:

a1) $d_1, d_2, ..., d_m$ belong to S;

a2) if a, b belong to S, then $\overline{ab}$ belongs to S too;

a3) only elements obtained by rules a1) and a2) applied

    a finite number of times belong to S.

Remark:

- there are m^k numbers of k digits in the sequence, for

  k = 1, 2, 3, ... ;

- to obtain from the k-digits number group the (k+1)-digits

  number group, just put first the digit $d_1$, second the digit $d_2$,

  ..., and the m-th time digit $d_m$ in the front of all k-digits

  numbers.

More general:  all digits $d_i$ can be replaced by numbers as
large as we want (therefore of many digits each), and also
m can be as large as we want.

Study these sequences.

**69. Square roots**:

0,1,1,1,2,2,2,2,2,3,3,3,3,3,3,3,4,4,4,4,4,4,4,4,4,5,5,5,5,

5,5,5,5,5,5,5,6,6,6,6,6,6,6,6,6,6,6,6,6,7,7,7,7,7,7,7,7,7,

7,7,7,7,7,7,8,8,8,8,8,8,8,8,8,8,8,8,8,8,8,8,8,9,9,9,9,9,9,

9,9,9,9,9,9,9,9,9,9,9,9,9,10,10,10,10,10,10,10,10,10,10,10,

10,10,10,10,10,10,10,10,10,...

(s$_q$(n) is the superior integer part of square root of n.)

Remark: this sequence is the natural sequence, where each number is repeated 2n+1 times,

because between n^2 (included) and (n+1)^2 (excluded) there are

(n+1)^2 - n^2 different numbers.

How many primes, squares, and cubes are in this sequence? What interesting properties has it?

## 70. Cubical roots:

0,1,1,1,1,1,1,1,2,2,2,2,2,2,2,2,2,2,2,2,2,2,2,2,2,2,2,3,3,3,

3,3,3,3,3,3,3,3,3,3,3,3,3,3,3,3,3,3,3,3,3,3,3,3,3,3,3,3,3,

3,3,3,3,4,4,4,4,4,4,4,4,4,4,4,4,4,4,4,4,4,4,4,4,4,4,4,4,4,

4,4,4,4,4,4,4,4,4,4,4,4,4,4,4,4,4,4,4,4,4,4,4,4,4,4,4,4,4,

4,4,4,4,4,...

(c$_q$(n) is the superior integer part of cubical root of n.)

Remark: this sequence is the natural sequence, where each number is repeated $3n^2 + 3n + 1$ times,

because between $n^3$ (included) and $(n+1)^3$ (excluded) there are

$(n+1)^3 - n^3$ different numbers.

## 71. m-power roots:

(m$_q$(n) is the superior integer part of m-power root of n.)

Remark: this sequence is the natural sequence, where each number is repeated $(n+1)^m - n^m$ times.

Study these sequences.

## 72. Numerical carpet:

has the general form

```
                .
                .
                .
                1
               1a1
              1aba1
             1abcba1
            1abcdcba1
           1abcdedcba1
          1abcdefedcba1
       ...1abcdefgfedcba1...
          1abcdefedcba1
           1abcdedcba1
            1abcdcba1
             1abcba1
              1aba1
               1a1
                1
                .
                .
                .
```

On the border of level 0, the elements are equal to "1";

   they form a rhomb.

Next, on the border of level 1, the elements are equal to "a",

   where "a" is the sum of all elements of the previous border;

   the "a"s form a rhomb too inside the previous one.

Next again, on the border of level 2, the elements are equal to "b",

   where "b" is the sum of all elements of the previous border;

   the "b"s form a rhomb too inside the previous one.

And so on...

The numerical carpet is symmetric and esthetic, in its middle g is the sum of all carpet numbers (the core).

Look at a few terms of the Numerical Carpet:

```
                              1

                              1
                             141
                              1

                              1
                           1  8  1
                        1  8 40  8  1
                           1  8  1
                              1

                              1
                           1  12   1
                        1  12 108 12  1
                     1 12 108 540 108 12 1
                        1  12 108 12  1
                           1  12   1
                              1

                              1
                           1  16   1
                        1  16 208  16  1
                    1  16 208 1872 208 16  1
                 1 16 208 1872 9360 1872 208 16 1
                    1 16  208 1872 208 16  1
                        1  16 208  16   1
                           1  16   1
                              1

                              1
                           1     20    1
                        1    20   340   20    1
                     1    20  340  4420  340   20   1
                  1   20  340 4420 39780 4420 340  20  1
              1 20 340 4420 39780 198900 39780 4420 340 20 1
                  1   20  340 4420 39780 4420 340  20   1
                     1   20  340  4420  340   20   1
                        1    20   340   20    1
                           1     20     1
                              1

                              .
                              .
                              .
```

Or, under other form:

```
1
1  4
1  8   40
1 12 108    504
1 16 208   1872    9360
1 20 340   4420   39780   198900
1 24 504   8568  111384  1002456    5012280
1 28 700  14700  249900  3248700   29238300  146191500
1 32 928  23200  487200  8282400  107671200  969040800  4845204000
..............................................................
    .
    .
    .
```

General Formula:

$$C(n,k) = 4n \prod_{i=1}^{k} (4n-4i+1) \text{ for } 1 \le k \le n,$$

and $C(n,0) = 1$.

Study this multi-sequence.

References:

a) Arizona State University, Hayden Library, "The Florentin Smarandache papers" special collection, Tempe, AZ 85287-1006, USA.

b) Student Conference, University of Craiova, Department of Mathematics, April 1979, "Some problems in number theory" by Florentin Smarandache.

c) F. Smarandache, "Collected Papers" (Vol. 1), Ed. Tempus, Bucharest, 1994 (to appear);

## 73. First Table:

6,10,14,18,26,30,38,42,42,54,62,74,74,90,...

(t(n) is the largest even number such that any other even number
not exceeding it is the sum of two of the first n odd primes.)

It helps to better understand Goldbach's conjecture:
 - if t(n) is unlimited, then the conjecture is true;
 - if t(n) is constant after a certain rank, then the
conjecture is false.

Also, the table gives how many times an even number is
written as a sum of two odd primes, and in what combinations.
Of course, $t(n) \le 2p_n$, where $p_n$ is the n-th odd prime,
n = 1, 2, 3, ... .

Here is the table:

```
       +    3   5   7  11  13  17  19  23  29  31  37  41  43  47
            ------------------------------------------- . . .
      3 |   6   8  10  14  16  20  22  26  32  34  40  44  46  50 .
      5 |      10  12  16  18  22  24  28  34  36  42  46  48  52 .
      7 |          14  18  20  24  26  30  36  38  44  48  50  54 .
     11 |              22  24  28  30  34  40  42  48  52  54  58 .
     13 |                  26  30  32  36  42  44  50  54  56  60 .
     17 |                      34  36  40  46  48  54  58  60  64 .
     19 |                          38  42  48  50  56  60  62  66 .
     23 |                              46  52  54  60  64  66  70 .
     29 |                                  58  60  66  70  72  76 .
     31 |                                      62  68  72  74  78 .
     37 |                                          74  78  80  84 .
     41 |                                              82  84  88 .
     43 |                                                  86  90 .
     47 |                                                      94 .
            .............................................
            .                                               .
            .                                                .
            .                                                 .
```

Study this table and table sequence.

## 74. Second table:

9,15,21,29,39,47,57,65,71,93,99,115,129,137,...

(v(n) is the largest odd number such that any odd number >= 9 not

exceeding it is the sum of three of the first n odd primes.)

It helps to better understand Goldbach's conjecture for

three primes:

- if v(n) is unlimited, then the conjecture is true;

- if v(n) is constant after a certain rank, then the

conjecture is false.

(Vinogradov proved in 1937 that any odd number greater than

3^(3^15) satisfies this conjecture.

But what about values less than 3^(3^15) ?)

Also, the table gives you in how many different combinations an odd number is written as a sum of three odd primes, and in what combinations.

Of course, v(n) <= 3$p_n$ , where $p_n$ is the n-th odd prime, n = 1,
2, 3, ... .  It is also generalized for the sum of m primes,

and how many times a number is written as a sum of m primes

(m > 2).

This is a 3-dimensional 14x14x14 table, that we can expose
only as 14 planar 14x14 tables (using the previous table):

```
       -----
      | 3   |
      |  +  |
      |     |   3   5   7  11  13  17  19  23  29  31  37  41  43  47
       -----  ----------------------------------------------------- . . .
         3  |  9  11  13  17  19  23  25  29  35  37  43  47  49  53 .
         5  |     13  15  19  21  25  27  31  37  39  45  49  51  55 .
         7  |         17  21  23  27  29  33  39  41  47  51  53  57 .
        11  |             25  27  31  33  37  43  45  51  55  57  61 .
        13  |                 29  33  35  39  45  47  53  57  59  63 .
        17  |                     37  39  43  49  51  57  61  63  67 .
        19  |                         41  45  51  53  59  63  65  69 .
        23  |                             49  55  57  63  67  69  73 .
        29  |                                 61  63  69  73  75  79 .
        31  |                                     65  71  75  77  81 .
        37  |                                         77  81  83  87 .
        41  |                                             85  87  91 .
        43  |                                                 89  93 .
        47  |                                                     97 .
             ............................................................
              .                                                       .
               .                                                       .
                .                                                       .

       -----
      | 5   |
      |  +  |
      |     |   3   5   7  11  13  17  19  23  29  31  37  41  43  47
       -----  ----------------------------------------------------- . . .
         3  | 11  13  15  19  21  25  27  31  37  39  45  49  51  55 .
         5  |     15  17  21  23  27  29  33  39  41  47  51  53  57 .
         7  |         19  23  25  29  31  35  41  43  49  53  55  59 .
        11  |             27  29  33  35  39  45  47  53  57  59  63 .
        13  |                 31  35  37  41  47  49  55  59  61  65 .
        17  |                     39  41  45  51  53  59  63  65  69 .
        19  |                         43  47  53  55  61  65  67  71 .
        23  |                             51  57  59  65  69  71  75 .
        29  |                                 63  65  71  75  77  81 .
        31  |                                     67  73  77  79  83 .
        37  |                                         79  83  85  89 .
        41  |                                             87  89  93 .
        43  |                                                 91  95 .
        47  |                                                     99 .
             ............................................................
              .                                                       .
               .                                                       .
                .                                                       .
```

```
      -----
     |  7  |
     |  +  |
     |     |   3   5   7  11  13  17  19  23  29  31  37  41  43  47
      -----  ---------------------------------------------------------  . . .
         3 | 13  15  17  21  23  27  29  33  39  41  47  51  53  57  .
         5 |     17  19  23  25  29  31  35  41  43  49  53  55  59  .
         7 |         21  25  27  31  33  37  43  45  51  55  57  61  .
        11 |             29  31  35  37  41  47  49  55  59  61  65  .
        13 |                 33  37  39  43  49  51  57  61  63  67  .
        17 |                     41  43  47  53  55  61  65  67  71  .
        19 |                         45  49  55  57  63  67  69  73  .
        23 |                             53  59  61  67  71  73  77  .
        29 |                                 65  67  73  77  79  83  .
        31 |                                     69  75  79  81  85  .
        37 |                                         81  85  87  91  .
        41 |                                             89  91  95  .
        43 |                                                 93  97  .
        47 |                                                    101  .
           ..............................................................
              .                                               .
              .                                                    .
              .                                                         .

      -----
     | 11  |
     |  +  |
     |     |   3   5   7  11  13  17  19  23  29  31  37  41  43  47
      -----  ---------------------------------------------------------  . . .
         3 | 17  19  21  25  27  31  33  37  43  45  51  55  57  61  .
         5 |     21  23  27  29  33  35  39  45  47  53  57  59  63  .
         7 |         25  29  31  35  37  41  47  49  55  59  61  65  .
        11 |             33  35  39  41  45  51  53  59  63  65  69  .
        13 |                 37  41  43  47  53  55  61  65  67  71  .
        17 |                     45  47  51  57  59  65  69  71  75  .
        19 |                         49  53  59  61  67  71  73  77  .
        23 |                             57  63  65  71  75  77  81  .
        29 |                                 69  71  77  81  83  87  .
        31 |                                     73  79  83  85  89  .
        37 |                                         85  89  91  95  .
        41 |                                             93  95  99  .
        43 |                                                 97 101  .
        47 |                                                    105  .
           ..............................................................
              .                                               .
              .                                                    .
              .                                                         .
```

```
      -----
     |13  |
     |  + |
     |    |    3    5    7   11   13   17   19   23   29   31   37   41   43   47
      ----- -----------------------------------------------
         3 |19   21   23   27   29   33   35   39   45   47   53   57   59   63  .
         5 |     23   25   29   31   35   37   41   47   49   55   59   61   65  .
         7 |          27   31   33   37   39   43   49   51   57   61   63   67  .
        11 |               35   37   41   43   47   53   55   61   65   67   71  .
        13 |                    39   43   45   49   55   57   63   67   69   73  .
        17 |                         47   49   53   59   61   67   71   73   77  .
        19 |                              51   55   61   63   69   73   75   79  .
        23 |                                   59   65   67   73   77   79   83  .
        29 |                                        71   73   79   83   85   89  .
        31 |                                             75   81   85   87   91  .
        37 |                                                  87   91   93   97  .
        41 |                                                       95   97  101  .
        43 |                                                            99  103  .
        47 |                                                                107  .
           ...............................................
           .                                                                      .
           .                                                                      .
           .                                                                      .

      -----
     |17  |
     |  + |
     |    |    3    5    7   11   13   17   19   23   29   31   37   41   43   47
      ----- -----------------------------------------------
         3 |23   25   27   31   33   37   39   43   49   51   57   61   63   67  .
         5 |     27   29   33   35   39   41   45   51   53   59   63   65   69  .
         7 |          31   35   37   41   43   47   53   55   61   65   67   71  .
        11 |               39   41   45   47   51   57   59   65   69   71   75  .
        13 |                    43   47   49   53   59   61   67   71   73   77  .
        17 |                         51   53   57   63   65   71   75   77   81  .
        19 |                              55   59   65   67   73   77   79   83  .
        23 |                                   63   69   71   77   81   83   87  .
        29 |                                        75   77   83   87   89   93  .
        31 |                                             79   85   89   91   95  .
        37 |                                                  91   95   97  101  .
        41 |                                                       99  101  105  .
        43 |                                                           103  107  .
        47 |                                                                111  .
           ...............................................
           .                                                                      .
           .                                                                      .
           .                                                                      .
```

```
 -----
|19  |
|  + |
|    |  |  3   5   7  11  13  17  19  23  29  31  37   41   43   47
 ----- ------------------------------------------------------------  . . .
    3 |25  27  29  33  35  39  41  45  51  53  59   63   65   69  .
    5 |    29  31  35  37  41  43  47  53  55  61   65   67   71  .
    7 |        33  37  39  43  45  49  55  57  63   67   69   73  .
   11 |            41  43  47  49  53  59  61  67   71   73   77  .
   13 |                45  49  51  55  61  63  69   73   75   79  .
   17 |                    53  55  59  65  67  73   77   79   83  .
   19 |                        57  61  67  69  75   79   81   85  .
   23 |                            65  71  73  79   83   85   89  .
   29 |                                77  79  85   89   91   95  .
   31 |                                    81  87   91   93   97  .
   37 |                                        93   97   99  103  .
   41 |                                            101  103  107  .
   43 |                                                 105  109  .
   47 |                                                      113  .
       ............................................................
         .                                                        .
          .                                                       .
           .                                                     .

 -----
|23  |
|  + |
|    |  |  3   5   7  11  13  17  19  23  29  31  37   41   43   47
 ----- ------------------------------------------------------------  . . .
    3 |29  31  33  37  39  43  45  49  55  57  63   67   69   73  .
    5 |    33  35  39  41  45  47  51  57  59  65   69   71   75  .
    7 |        37  41  43  47  49  53  59  61  67   71   73   77  .
   11 |            45  47  51  53  57  63  65  71   75   77   81  .
   13 |                49  53  55  59  65  67  73   77   79   83  .
   17 |                    57  59  63  69  71  77   81   83   87  .
   19 |                        61  65  71  73  79   83   85   89  .
   23 |                            69  75  77  83   87   89   93  .
   29 |                                81  83  89   93   95   99  .
   31 |                                    85  91   95   97  101  .
   37 |                                        97  101  103  107  .
   41 |                                            105  107  111  .
   43 |                                                 109  113  .
   47 |                                                      117  .
       ............................................................
         .                                                        .
          .                                                       .
           .                                                     .
```

```
       -----
      |29   |
      |  +  |
      |     |  3   5   7  11  13  17  19  23  29  31   37   41   43   47
       ----- ---------------------------------------------------------------  . . .
         3  |35  37  39  43  45  49  51  55  61  63   69   73   75   79   .
         5  |    39  41  45  47  51  53  57  63  65   71   75   77   81   .
         7  |        43  47  49  53  55  59  65  67   73   77   79   83   .
        11  |            51  53  57  59  63  69  71   77   81   83   87   .
        13  |                55  59  61  65  71  73   79   83   85   89   .
        17  |                    63  65  69  75  77   83   87   89   93   .
        19  |                        67  71  77  79   85   89   91   95   .
        23  |                            75  81  83   89   93   95   99   .
        29  |                                87  89   95   99  101  105   .
        31  |                                    91   97  101  103  107   .
        37  |                                        103  107  109  113   .
        41  |                                             111  113  117   .
        43  |                                                  115  119   .
        47  |                                                       123   .
            .................................................................
             .                                                        .
              .                                                         .
               .                                                          .

       -----
      |31   |
      |  +  |
      |     |  3   5   7  11  13  17  19  23  29  31   37   41   43   47
       ----- ---------------------------------------------------------------  . . .
         3  |37  39  41  45  47  51  53  57  63  65   71   75   77   81   .
         5  |    41  43  47  49  53  55  59  65  67   73   77   79   83   .
         7  |        45  49  51  55  57  61  67  69   75   79   81   85   .
        11  |            53  55  59  61  65  71  73   79   83   85   89   .
        13  |                57  61  63  67  73  75   81   85   87   91   .
        17  |                    65  67  71  77  79   85   89   91   95   .
        19  |                        69  73  79  81   87   91   93   97   .
        23  |                            77  83  85   91   95   97  101   .
        29  |                                89  91   97  101  103  107   .
        31  |                                    93   99  103  105  109   .
        37  |                                        105  109  111  115   .
        41  |                                             113  115  119   .
        43  |                                                  117  121   .
        47  |                                                       125   .
            .................................................................
             .                                                        .
              .                                                         .
               .                                                          .
```

```
     -----
    |37  |
    |  + |
    |    |   3   5   7  11  13  17  19  23  29  31  37  41  43  47
     ----- ------------------------------------------------  . . .
       3 |43  45  47  51  53  57  59  63  69  71  77  81  83  87  .
       5 |    47  49  53  55  59  61  65  71  73  79  83  85  89  .
       7 |        51  55  57  61  63  67  73  75  81  85  87  91  .
      11 |            59  61  65  67  71  77  79  85  89  91  95  .
      13 |                63  67  69  73  79  81  87  91  93  97  .
      17 |                    71  73  77  83  85  91  95  97 101  .
      19 |                        75  79  85  87  93  97  99 103  .
      23 |                            83  89  91  97 101 103 107  .
      29 |                                95  97 103 107 109 113  .
      31 |                                    99 105 109 111 115  .
      37 |                                       111 115 117 121  .
      41 |                                           119 121 125  .
      43 |                                               123 127  .
      47 |                                                   131  .
          .................................................
          .                                               .
          .                                                 .
          .                                                   .

     -----
    |41  |
    |  + |
    |    |   3   5   7  11  13  17  19  23  29  31  37  41  43  47
     ----- ------------------------------------------------  . . .
       3 |47  49  51  55  57  61  63  67  73  75  81  85  87  91  .
       5 |    51  53  57  59  63  65  69  75  77  83  87  89  93  .
       7 |        55  59  61  65  67  71  77  79  85  89  91  95  .
      11 |            63  65  69  71  75  81  83  89  93  95  99  .
      13 |                67  71  73  77  83  85  91  95  97 101  .
      17 |                    75  77  81  87  89  95  99 101 105  .
      19 |                        79  83  89  91  97 101 103 107  .
      23 |                            87  93  95 101 105 107 111  .
      29 |                                99 101 107 111 113 117  .
      31 |                                   103 109 113 115 119  .
      37 |                                       115 119 121 125  .
      41 |                                           123 125 129  .
      43 |                                               127 131  .
      47 |                                                   135  .
          .................................................
          .                                               .
          .                                                 .
          .                                                   .
```

```
       -----
      |43   |
      |  +  |
      |     |   3   5   7  11  13  17  19  23   29   31   37   41   43   47
       ----- ------------------------------------------------------------   . . .
          3 |49  51  53  57  59  63  65  69   75   77   83   87   89   93  .
          5 |    53  55  59  61  65  67  71   77   79   85   89   91   95  .
          7 |        57  61  63  67  69  73   79   81   87   91   93   97  .
         11 |            65  67  71  73  77   83   85   91   95   97  101  .
         13 |                69  73  75  79   85   87   93   97   99  103  .
         17 |                    77  79  83   89   91   97  101  103  107  .
         19 |                        81  85   91   93   99  103  105  109  .
         23 |                            89   95   97  103  107  109  113  .
         29 |                                101  103  109  113  115  119  .
         31 |                                     105  111  115  117  121  .
         37 |                                          117  121  123  127  .
         41 |                                               125  127  131  .
         43 |                                                    129  133  .
         47 |                                                         137  .
            ...................................................
             .                                                     .
             .                                                     .
             .                                                       .

       -----
      |47   |
      |  +  |
      |     |   3   5   7  11  13  17  19  23   29   31   37   41   43   47
       ----- ------------------------------------------------------------   . . .
          3 |53  55  57  61  63  67  69  73   79   81   87   91   93   97  .
          5 |    57  59  63  65  69  71  75   81   83   89   93   95   99  .
          7 |        61  65  67  71  73  77   83   85   91   95   97  101  .
         11 |            69  71  75  77  81   87   89   95   99  101  105  .
         13 |                73  77  79  83   89   91   97  101  103  107  .
         17 |                    81  83  87   93   95  101  105  107  111  .
         19 |                        85  89   95   97  103  107  109  113  .
         23 |                            93   99  101  107  111  113  117  .
         29 |                                105  107  113  117  119  123  .
         31 |                                     109  115  119  121  125  .
         37 |                                          121  125  127  131  .
         41 |                                               129  131  135  .
         43 |                                                    133  137  .
         47 |                                                         141  .
            ...................................................
             .                                                     .
             .                                                     .
             .                                                       .
```

Second table sequence:

0,0,0,0,1,2,4,4,6,7,9,10,11,15,17,16,19,19,23,25,26,26,28,

33,32,35,43,39,40,43,43,...

(a(2k+1) represents the number of different combinations

such that 2k+1 is written as a sum of three odd primes.)

This sequence is deduced from the second table.

Study the second table and the second table sequence.

References:

a) Florentin Smarandache, "Problems with and without ... problems!",

   Ed. Somipress, Fes, Morocco, 1983;

b) Arizona State University, Hayden Library, "The Florentin

   Smarandache papers" special collection, Tempe, AZ 85287-

   1006, USA.